\documentclass[12pt]{amsart}

\usepackage{eucal,url,amssymb,stmaryrd,
enumerate,amscd,
}
\usepackage[pagebackref]{hyperref}  

\numberwithin{equation}{section}

\newtheorem{thm}{Theorem}[section]

\newtheorem{lem}[thm]{Lemma}

\newtheorem{prop}[thm]{Proposition}

\newtheorem{dfn}[thm]{Definition}

\newtheorem{exm}[thm]{Example}

\newcommand{\m}{\frak{m}}

\begin{document}

\title{When is tight closure determined by the test ideal?}
\author{Janet C. Vassilev}
\address{The University of California, Riverside, CA 92521}
\email{jvassil@ucr.edu}
\author{Adela N. Vraciu}
\address{University of South Carolina, Columbia, SC 29205}
\email{vraciu@math.sc.edu} \keywords{ basically-full, tight closure,
test ideal}

\begin{abstract}{We characterize the rings in which
the equality $(\tau I:\tau)= I^*$ holds for every ideal $I \subset R$. Under certain assumptions,
these rings must be either weakly F-regular or one-dimensional.}\end{abstract} \maketitle

\section{Introduction}
Test ideals play a major role in the theory of tight closure. The tight closure of arbitrary ideals is very
difficult to compute, even in relatively simple rings, but the test ideal can be frequently computed,
especially in Gorenstein rings. Moreover, test ideals encode geometric information about the nature of
the singularity of the ring.
We recall the definitions and basic facts.

Throughout this paper, $(R, \frak{m})$ is a local domain of characteristic $p$.  We
denote positive integer powers of $p$ by $q$.
\begin{dfn} Let $I\subset R$ be an ideal. We say that $x \in R$
is in the { \it tight closure}, $I^*$, of $I$ if there is a $c \neq
0$ such that $cx^q \in I^{[q]}=(\{i^q|i \in I\})$.  We say that $I$ is
{\it tightly closed} if $I=I^*$.
\end{dfn}
\begin{dfn}
The {\it test ideal} $\tau $ is defined by
$$
\tau = \bigcap _{I \subset R} (I:I^*),
$$ where $I$ runs over all ideals $I \subset R$.
\end{dfn}
The fact that $\tau \ne (0)$ is a highly nontrivial and important result.

It is clear from the definition that $I^* \subseteq I : \tau $, and thus $I : \tau $ provides an upper bound for
tight closure. This bound can be somewhat refined with additional assumptions on the ring, as shown in the
following result of the second author:
\begin{thm}\cite{Vr} If $R$ is a complete domain of characteristic $p$, then  the test
ideal is a strong test ideal, i.e. we have $\tau I = \tau I^*$, and thus $I^* \subseteq (\tau I : \tau$) for all
ideals $I \subset R$.
\end{thm}
Also in the case when the test ideal is the maximal ideal, a theorem of Hara and Smith \cite{HS} says that over
any local ring with test ideal equal to the maximal ideal, the test ideal is the strong test ideal.

We say that tight closure is determined by the test ideal if the equality $I^* = \tau I : \tau $ holds.
This is known to hold if the ideal $I$ is generated by a system of parameters in a Gorenstein ring $R$
(Corollary 4.2 (2) in \cite{Hu1}), or, more generally, if $I$ is an ideal of finite projective dimension
in a Gorenstein ring (Theorem 1(a) in \cite{Vlink}). The main result of this paper shows, under certain assumptions,
that the equality {\it cannot}
hold for all ideals $I \subset R$ unless the ring is weakly F-regular (i.e. $\tau = R$) or one-dimensional.

In a similar vein, we mention a result of Yao [Thm. 2.5 (ii)] in~\cite{Ya}, which states that if $R$ has finite
Frobenius representation type, then there exists a finitely generated $R$-module $N$ such that
$I^* = (IN:_R N)$ for every ideal $I\subset R$. Thus, our result indicates that the $R$-module $N$ cannot be an
ideal unless $R$ is one-dimensional or weakly F-regular.

A different motivation for our work comes from the following result of Heinzer, Ratliff and Rush
in \cite{HRR}[Theorem 7.5]:
\begin{thm}\label{HHR}  Let $(R,\frak{m})$ be a local ring.  A necessary and
sufficient condition for every nonzero $\frak{m}$-primary ideal of
$R$ to be basically full is that $\frak{m}$ is principal and that
$R$ is a principal ideal ring.
\end{thm}
In their terminology, an $\frak{m}$-primary ideal $I$ is { \it
basically full} if $(\frak{m}I:\frak{m})=I$, and for any ideal $I$,
$(\frak{m}I:\frak{m})$ is called the {\it basically full closure} of
$I$. The original motivation for this paper was the desire to extend
Theorem~\ref{HHR} to take tight closure into account. Thus, we asked
the question: when is $(\frak{m}I:\frak{m})=I^*$ for all
$\frak{m}$-primary ideals $I\subset R$? Theorem~\ref{thm1} shows
that (under certain assumptions) this is the case if and only if $R$
is one-dimensional.

\section{$*$-$T$-basically full ideals}

We extend the definition of the basically full closure in \cite{HRR}
using any ideal $T$ to define the {\it $T$-basically full closure}
of an ideal $I$ to be $I^{T bf}=(T I:T)$. This is a true closure
operation as:

\begin{prop}  Let $(R,\frak{m})$ be a local domain. The $T$-basically full closure of an ideal $(T I:T)$ is a
closure operation satisfying:
\begin{enumerate}
\item $I \subseteq I^{T bf}$.

\item If $I \subseteq J$, then $I^{T bf} \subseteq J^{T bf}$.

\item $(I^{T bf})^{T bf}=I^{T bf}$.

\item $I^{T bf}J^{T bf} \subseteq (IJ)^{T bf}$.

\end{enumerate}
\end{prop}

{\bf Proof:}  (1) and (2) are clear.

For (3) note that if $I$ is any ideal and $I^{T bf}=(T I:T)$, $T
I^{T bf}=T (TI:T) = T I$, hence $(T I^{T bf}: T) =(T I:T)=I^{T bf}$.

For (4), note that $$I^{T bf}J^{T bf} =(TI:T)(T J:T) \subseteq (T (T
I:T)(T J:T):T)$$
$$\quad \quad \quad\subseteq
 (T I(TJ:T):T) \subseteq (T IJ:T)=(IJ)^{T bf}.$$  \qed

We want to determine the domains which satisfy $I^{T bf}=I^{*}$ for
all $\m$-primary ideals $I$.  This prompts the following definition:

\begin{dfn} $I$ is $T$-basically full if $I^{T bf}=I$.  We
will say that $I$ is $*$-$T$-basically full if $I^{T bf}=I^*$.
\end{dfn}

\begin{thm} \label{*tau}  Let $(R,\frak{m})$ be a complete local normal Cohen Macaulay domain of positive
characteristic with perfect residue field having a canonical module
and let $\tau$ be the test ideal. If $T$ is an ideal of grade at
least two, then every $\m$-primary ideal is $*$-$T$-basically full
if and only if $R$ is weakly $F$-regular and $T=R$. In particular,
every $\m$-primary ideal is $*$-$\tau$-basically full if and only if
$R$ is weakly F-regular.
\end{thm}

Before the proof, note that the normal assumption is necessary.   If $R$ is a one-dimensional domain (in which case normal is equivalent to regular and therefore it is also equivalent to weakly F-regular) the following Proposition shows that every
$\frak{m}$-primary ideal is $*$-$\tau$-basically full.

\begin{prop} \label{1dim} Assume that $(R,\frak{m})$ is complete domain with infinite residue field. If $R$ has Krull
dimension one, then every $\frak{m}$-primary ideal is
$*$-$\tau$-basically full.
\end{prop}

{\bf Proof:}  In a one-dimensional domain with infinite residue field, every $\frak{m}$-primary ideal $I$ has a principal minimal reduction $(x)$. For principal ideals, tight closure is the same as integral closure, and it follows that $I^*=(x)^*$
(see \cite{Hu}[Example 1.6.2]). 

We have $$(\tau I:\tau) \subseteq (\tau
I^*:\tau)=(\tau(x)^*:\tau) = (\tau (x):\tau)\subseteq\overline{(x)}=(x)^*=I^*.$$
The equality $\tau(x)^*:\tau=\tau(x):\tau$ uses the fact that $\tau$ is a strong test ideal, and the inclusion $\tau (x):\tau \subseteq \overline{(x)}$ uses the determinant trick.

 Since
$\tau$ is a strong test ideal, the inclusion  $I^* \subseteq (\tau I:\tau)$ also holds. 
\medskip

Note that in the case of a one-dimensional domain, the only non
$\m$-primary ideals are $(0)$ and $R$, and thus Theorem~\ref{*tau}
shows that in this case $I^*=\tau I :\tau $ holds for every ideal
$I$, hence the tight closure of every ideal in a one-dimensional
domain is determined by the test ideal.

\medskip

{\bf Proof of \ref{*tau}:}
Note that the last statement follows from the previous one, since in an excellent normal ring the test ideal always has
depth at least two (see Theorem 6.2 in \cite{HH2}).

 One implication is clear: if $R$ is
weakly $F$-regular and $T=R$, then $I^* = I =TI:T$ for every ideal $I \subset R$.

Conversely, suppose that all $\frak{m}$-primary ideals are
$*$-$T$-basically full.
It is enough to prove that $T$ must be a principal ideal, because then the grade assumption implies that $T=R$, and thus
$I^*=TI:T=I$ for every $\frak{m}$-primary ideal $I$, which implies that $R$ is weakly F-regular.

Assume by contradiction that the minimal number of
generators of $T$ is $\nu (T)=n\ge 2$, and write
$T=(y_1,y_2,\ldots,y_n)$.

Following the argument (2.2.1), the proof of Theorem 2.2.2 in \cite{Ab},
suppose $J$ is a canonical ideal and choose $x_1, \ldots, x_d$ a system of
parameters for $R$ such that $x_1 \in J$, and $x_2, \ldots, x_d$ form a
regular sequence modulo $J$. Note that $\frak{a}=(J,x_2,\ldots x_d)$ is an
irreducible ideal, and $\frak{a}_t=(x_1^{t-1}J,x_2^t,\ldots x_d^t)$ are
irreducible for all $t\ge 2$.  Let $v$ denote the socle element of
$\frak{a}$, i.e. $(\frak{a}:\frak{m})=(\frak{a},v)$ and let  $v_t=(x_1\cdots
x_d)^{t-1}v$ be the socle element of $\frak{a}_t$.

By Matlis duality, we have
$$\lambda \left(\frac{(\frak{a}_t :\frak{m}T)}{(\frak{a}_t :T)}\right)=
\lambda\left(\frac{\frak{a}_t + T}{\frak{a}_t + \frak{m}T}\right)=
\lambda \left( \frac{T}{\frak{a}_t \cap T+ \frak{m}T}\right).$$
Note that $\frak{a}_t \cap T\subseteq \frak{m}^t \cap T\subseteq \m T$ for $t\gg 0$ by the
Artin-Rees Lemma, and therefore this length is equal to one if and only if $T$ is a principal ideal.

Fix a $t_0$ large enough so that
$$\lambda \left(\frac{(\frak{a}_{t_0} :\frak{m}T)}{(\frak{a}_{t_0} :T)}\right)\ge 2
$$
 and choose $u_1, u_2 \in (\frak{a}_{t_0} :\frak{m}T)$
such that their images are linearly independent in the vector space
$(\frak{a}_{t_0} :\frak{m}T)/(\frak{a}_{t_0} :T)$.

 Note that for all $t \ge 1$, $(x_1 \cdots x_d)^{t} u_1, (x_1
\cdots x_d)^{t} u_2\in (\frak{a}_{t_0+t}: \frak{m} T)$, and their
images in  $(\frak{a}_{t_0+t}: \frak{m} T)/ (\frak{a}_{t_0 + t}: T)
$ are linearly independent, because the map $R/\frak{a}_{t_0}
\rightarrow R/\frak{a}_{t_0 + t}$ given by multiplication by $(x_1
\cdots x_d)^{t}$ is injective.

Consider the ideals $$ I_{t1}=(\frak{a}_{t_0 + t}, (x_1 \cdots
x_d)^{t} u_2), \ I_{t2}=(\frak{a}_{t_0 + t},(x_1 \cdots x_d)^{t}
u_1).$$

We claim that $(x_1 \cdots x_d)^{t} u_i\in TI_{ti} : T$ for
$i=1, 2$ when $t\gg 0$.

The key point in the proof of the claim is the observation that the
assumption that $T$ has grade at least two implies that we can choose
$x_1, \ldots, x_d$  so that at least two of them belong to $T$ (by prime
avoidance). With the $x$'s chosen this way, we have
\begin{equation}\label{Teq} (x_1\ldots x_d)^{t}\frak{a}_{t_0}\subseteq
T\frak{a}_{t_0+t}.\end{equation}

For $j=1, \ldots, n$, we have $ u_i y_j \in \frak{a}_{t_0} : \frak{m}
=(\frak{a}_{t_0}, v_{t_0})$. Moreover, since $ u_i\notin \frak{a}_{t_0}
:T$, for each $i$ there exists a $j=j_i$ such that $ u_i y_{j_i} \notin
\frak{a}_{t_0}$, so that we can write $ u_i y_{j_i}= \alpha v_{t_0} \
(\mathrm{mod} \ \frak{a}_{t_0})$, where $\alpha $ is a unit. This shows
that $v_{t_0}\in (\frak{a}_{t_0}, Tu_i)$ for all $i$, and therefore
$$(x_1\cdots x_d)^t v_{t_0} \in ((x_1\ldots
x_d)^{t}\frak{a}_{t_0},T(x_1\ldots x_d)^{t} u_i) \subseteq
T(\frak{a}_{t_0+t},(x_1\ldots x_d)^{t} u_i)$$ when $t \gg 0$ by Equation
\ref{Teq}. Combining Equation \ref{Teq} with the above chain of
containments, we conclude that $(x_1\cdots x_d)^{t} (\frak{a}_{t_0},
v_{t_0}) \subseteq TI_{ti}$ for each $i$ (when $t \gg 0$), which finishes
the proof of the claim, since $(x_1 \cdots x_d)^{t} u_i T \subseteq
(x_1\cdots x_d)^{t} (\frak{a}_{t_0}, v_{t_0}) $ by the choice of $u_1,
u_2$.

Since the ideals $I_{ti}$ are assumed to be $*-T-$basically full, we have
$$(x_1 \cdots x_d)^{t} u_i\in (TI_{ti}:T)=I_{ti}^*.$$ Note that the same argument works
when $u_1$ is replaced by $u_1 + \alpha u_2$, where $\alpha \in R$ is
arbitrary. Therefore,
 Lemma~\ref{special} can be applied to see that $(x_1 \cdots x_d)^{t} u_i \in \frak{a}_{t_0 +t}^*$.
 This is a contradiction, since $(x_1 \cdots x_d)^{t} u_i\notin \frak{a}_{t_0 +t} :
 T$ and
 $T \frak{a}_{t_0 +t}: T \subsetneq \frak{a}_{t_0 +t} : T$ which implies
 that
 $(x_1 \cdots x_d)^{t} u_i\notin \frak{a}_{t_0 +t}^*=T \frak{a}_{t_0 +t}: T$.
\qed

\begin{lem}\label{special}
Assume that $R$ is a complete normal domain of positive
characteristic $p$, with perfect residue field. If $I \subset R$ is
an ideal, and $f, g \in R$ are such that $f \in (I, g)^*$ and $g \in
(I, f+\alpha g)^*$ for all $\alpha \in R$, then $f, g \in I^*$.
\end{lem}
{\bf Proof:} Theorem 2.1 in \cite{HV} shows that $(I, g)^*=(I, g) + (I,
g)^{*sp}$, and therefore there exists an $\alpha \in R$, a $q_0=p^{e_0}$,
and a $c \in R^o$ such that $c(f+ \alpha g)^q =bg^q \ \mathrm{mod}
I^{[q]}$ for all $q=p^e$, with $b \in \m^{q/q_0}$. On the other hand,
there exists $c' \in R^o$ such that $c'g^q =d(f+\alpha g)^q \ \mathrm{mod}
I^{[q]}$. Combining these two equations, we get $cc' (f+\alpha g)^q =bd(f
+ \alpha g)^q \ \mathrm{mod} I^{[q]}$. Since $bd \in \m^{q/q_0}$ and $c,
c'$ are fixed, Proposition 2.4 in \cite{Ab1} shows that $f+ \alpha g \in
I^*$. Since $g \in (I, f+\alpha g)^*$, we also get $g \in I^*$, and since
$f \in (I, g)^*$ we now get $f \in I^*$ as well. \qed

We cannot remove the assumption of perfect residue field in Lemma
\ref{special}.  Consider the following example motivated by \cite{Ep}[p.
381]:

\begin{exm}\label{neil}  Let
$R=\displaystyle\frac{\mathbb{Z}/p\mathbb{Z}(u,v,w)[[x,y,z]]}{(ux^p+vy^p+wz^p)}$
which is a $2$-dimensional, Gorenstein normal domain as remarked by
Epstein in \cite{Ep}. Let $I=(x^2,y^2,z)$.  $x \in (I,y)^F \subseteq
(I,y)^*$ and for all $a \in R$, $y \in (I,x+ay)^F \subseteq (I,x+ay)^*$.
However, $x,y \notin I^*=(xy,x^2,y^2,z)$.  Hence Lemma \ref{special}
requires a perfect residue field.
\end{exm}

It may be however that Theorem \ref{*tau} holds when the residue field is
not perfect, as the above ring does not satisfy $I^*= (\tau I:\tau)$
for all $I$.  To see this we will compute the test ideal for $R$ and
exhibit the offending $\frak{m}$-primary ideal $I$.

We claim that for all $t \ge p$ we have $(y^t, z^t)^*=(y^t,
z^t)+\frak{m}^{2t-1} =(y^t, z^t):\frak{m}^{p-1}$, and thus $\tau
=\frak{m}^{p-1}$.

In order to prove the first equality, it is enough to consider monomials of the form $x^ky^rz^s$, with $k\le p-1$.
Note that $x^ky^rz^s\in (y^t, z^t) \Leftrightarrow x^k \in (y^{t-r}, z^{t-s})^*\Leftrightarrow x^{kp}
\in (y^{(t-r)p}, z^{(t-s)p})$.

 We have $$x^{kp}=(-\frac{1}{u}(v y^p+w
z^p))^k=\frac{(-1)^k}{u^k} \underset{i=0}{\overset{k}{\sum}}\ \left(\begin{array}{c}  k \\ i\\ \end{array}\right)
v^iw^{k-i}y^{ip}z^{(k-i)p},$$
and thus the tight closure membership can be tested inside the regular ring $k[[y, z]]$. Since $k \le p-1$,
none of the binomial coefficients $\displaystyle \left(\begin{array}{c}  k \\ i\\ \end{array}\right)$ is equal to
zero, and thus we see that $x^k \in (y^{t-s}, z^{t-s})^* \Leftrightarrow $ for all $i =0, \ldots, k$ we have either
$i \ge t-r$, or $k-i \ge t-s$. This amounts to $ k \ge 2t -r -s -1$, which proves the first equality.

For the second equality, it is enough to show that $(y^t, z^t):\frak{m}^{2t-1}=(y^t, z^t)+\m^{p-1}$
(since the ring is Gorenstein). It is easy to see that $(y^t, z^t):\frak{m}=
(x^t, y^t, y^{t-1}z^{t-1}x^{p-1})=(y^t, z^t)+\frak{m}^{2t+p-3}$, and one can check by induction on $l$ that
$(y^t, z^t):\frak{m}^l=(y^t, z^t)+\frak{m}^{2t+p-2-l}$ for all $l\ge 1$. Taking $l=2t-1$ yields the desired conclusion.

Now we show that for $p\ge 5$, these rings do not have the property that every $\m$-primary ideal is
$*$-$\tau$-basically full.

Consider $I=(x^3,y^3,z^3)$ and $r=xy^2$. We have $r \in
(\frak{m}^{p-1} I:\frak{m}^{p-1})$; however, $r \notin I^*$. Using
the relation, $z^p=-\frac{u}{w} x^p-\frac{v}{w}y^p$, we see that
$$z^{3p}=-\frac{u^3}{w^3}x^{3p}-3\frac{u^2v}{w^3}x^{2p}y^p-3\frac{uv^2}{w^3}x^py^{2p}-\frac{v^3}{w^3}y^{3p}.$$
Hence if $xy^2 \in I^*$, then $x^py^{2p} \in
(x^{3p},y^{3p},z^{3p})^*=(x^{3p},y^{3p},u^2vx^{2p}y^p+uv^2x^{p}y^{2p})^*$,
implying that $x^py^{2p} \in (x^{3p},y^{3p},x^{2p}y^p)^*$ in the
regular ring $k[[x, y]]$.  This leads to a contradiction.

For $p=3$, consider the ideal $I=(x^4,y^4,z^4)$ and $r=xy^3$.  Note
that $r \in (\frak{m}^2 I:\frak{m}^2)$ if and only if $x^3y^3 \in I$
which is the case as
$x^3y^3=\frac{u}{v}x^6+\frac{v}{u}y^6+2\frac{w^2}{uv}z^6 \in
(x^4,y^4,z^4)$.  As in the argument above, $xy^3 \in I^*$ is
equivalent to $xy^3 \in (x^4,y^4,x^3y)^*$ in the regular ring
$k[[x,y]]$.  Again this leads to a contradiction.

For $p=2$, we have not found an ideal $I$ for which
$(\frak{m}I:\frak{m}) \neq I^*$.  For characteristic 3 and higher,
the ring $k[[x,y,z]]/(x^2+y^2+z^2)$ is $F$-rational. Certainly,
$\displaystyle\frac{\mathbb{Z}/2\mathbb{Z}(u,v,w)[[x,y,z]]}{(ux^2+vy^2+wz^2)}$
is not as $(y,z)^*=\frak{m}$, but it may be that this ring satisfies
$(\frak{m}I:\frak{m})=I^*$ for all $\frak{m}$-primary ideals $I$.

Since we do not know whether the conclusion of Theorem~\ref{*tau}
holds in the absence of the perfect residue field assumption, it is
worth pointing out that if $R$ is a Gorenstein ring and $T$ is an ideal of grade at least two
such that for every $\m$-primary ideal $I$ we have $I^*=TI:T$, then
$T$ is forced to be the test ideal. In particular, this shows that
the full conclusion of Theorem~\ref{*tau} holds for the rings
considered in Example~\ref{neil} for $p \ge 3$.

Before the following proposition, we recall a definition of Hochster in
his paper on Cyclic Purity \cite{Ho}[Definition 1.1]:

\begin{dfn}
A local Noetherian ring $(R,\m)$ is approximately Gorenstein if for every
$N >0$, there is an $\m$-primary ideal $I \subseteq \m^N$ which is
irreducible.
\end{dfn}

Hochster noted in \cite{Ho}[Remark 4.8(b)] that a generically
Gorenstein Cohen Macaulay ring with a canonical module is
approximately Gorenstein.

\begin{prop}
Assume that $(R, \frak{m})$ is an approximately Gorenstein ring with test ideal $\tau $. If $T$
is an ideal such that every $\frak{m}$-primary ideal of $R$ is
$*$-$T$-basically full, then $T\subseteq \tau$.

If we moreover assume that $R$ is Gorenstein and $T$ has grade at least two, then
$T=\tau$.
\end{prop}
Note that the assumption that $T$ has grade at least two is
necessary in the second part of the Proposition. If $R$ is a weakly
F-regular domain and $T$ is a principal ideal, then $I^*=I=TI:T$ for
every $I$, but $T\ne \tau =R$.
\medskip

{\bf Proof:} Let $\{ \frak{a}_t\}$ be a sequence of irreducible ideals cofinal with
the powers of $\frak{m}$. We have
$$
\tau = \bigcap _t \frak{a}_t : \frak{a}_t^*=
\bigcap _t \frak{a}_t : (T\frak{a}_t :T)\supseteq \bigcap _t \frak{a}_t : (\frak{a}_t :T) =
\bigcap _t (\frak{a}_t + T)=T.
$$

Now assume that $R$ is Gorenstein and $T$ has grade at least two. Let $x_1, \ldots, x_d$ be a system of parameters for $R$ such that at least two of them belong to $T$.

Let $\frak{a}_t = (x_1^t, x_2^t, \ldots, x_d^t)$. We claim that $\frak{a}_t^* = \frak{a}_t :T$
for all $t$. Assuming the claim, we obtain
$$
T=\bigcap_t (\frak{a}_t + T) = \bigcap_t (\frak{a}_t: (\frak{a}_t : T))=\bigcap _t (\frak{a}_t : \frak{a}_t ^*) = \tau.
$$

In order to prove the claim, consider $u \in \frak{a}_t: T$. Since at
least two out of $x_1, \ldots, x_d$ belong to $T$, we have $(x_1\cdots
x_d)\frak{a}_t \subseteq T\frak{a}_{t+1}$, and therefore $(x_1\cdots
x_d)u\in T\frak{a}_{t+1} : T = \frak{a}_{t+1}^*$. Thus, $u\in  
\frak{a}_{t+1}^*:(x_1\cdots x_d) \subseteq \frak{a}_t^*$. \qed

\section{When the basically full closure and the tight closure correspond for all $m$-primary ideals}

In this section we extend the definition of basically full closure of \cite{HRR} in a slightly different direction, using the maximal ideal $\m$ instead of $\tau $ (thus, this version is closer to the original one in \cite{HRR}).

\begin{dfn} We will say an $\frak{m}$-primary ideal $I$
is {\it $*$-basically full} if $(\frak{m}I:\frak{m})=I^*$.
\end{dfn}

\begin{thm}\label{thm1} Let $(R,\frak{m})$ be a local Cohen Macaulay domain.

(a). If $R$ is a one-dimensional ring with
test ideal equal to $\frak{m}$, then all $\frak{m}$ primary
ideals are $*$-basically full.

(b). Assume in addition that $R$ is normal and has perfect residue field.
Then $R$ is a one-dimensional ring if and only if all $\frak{m}$ primary ideals are $*$-basically full.
\end{thm}
{\bf Proof:} (a). Follows from the same proof as in Proposition
\ref{1dim} (using the fact that when $\m$ is the test ideal, it is a
strong test ideal even without assuming that the ring is complete).

(b). Assume that all $\frak{m}$ primary ideals are $*$-basically full. Note that this assumption implies that $\tau =\m$ or $\tau =R$, since
$$\tau = \bigcap_{I \subset R} (I: I^*)=\bigcap_{I \subset R} (I:(\m I : \m)) \supseteq \bigcap _{I \subset R} (I : (I : \m))\supseteq \m.$$
If $\tau =\m$, then we are under the assumptions of
Theorem~\ref{1dim}, and thus $\tau =R$, which is a contradiction. If
$\tau =R$, then we have $(\m I : \m)=I^* =I$, i.e. every
$\m$-primary ideal is basically full, and Theorem~\ref{HHR} applies
to show that $R$ must have dimension one. \qed

To see some examples we will use the following theorem from the
first author's thesis \cite{Va}:

\begin{thm}  Let $(R,\frak{m})$ be a one-dimensional domain.  The test
ideal of $R$ is equal to the conductor, i.e. $\tau=\frak{c}=\{c \in
R|\phi(1)=c,\phi \in \text{Hom}_R(\overline{R},R)\}$.
\end{thm}

Note, in a one-dimensional local semigroup ring, the semigroup is a
sub-semigroup of $N_0$.  For each sub-semigroup $S$ of $N_0$, there
is a smallest $m$ such that for all $i \geq m$, $i \in S$. The
conductor of such a one dimensional semigroup ring is
$\frak{c}=<t^m, t^{m+1},t^{m+2}, \ldots>$, \cite{Ei}[Exercise
21.11].

\begin{exm} The rings $k[[t^2,t^3]]$ and $k[[t^3,t^4,t^5]]$ are one dimensional domains with
test ideal equal to the maximal ideal.  Every $\frak{m}$-primary
ideal of each ring is $*$-basically full. \end{exm}

\begin{exm}  The test ideal of $k[[t^2,t^5]]$ is $(t^4,t^5)$ , hence there are ideals in $k[[t^2,t^5]]$ which
are not $*$-basically full.  For example,
$(t^4)^*=(t^4):(t^4,t^5)=(t^4,t^5)$ and \newline
$(\frak{m}(t^4):\frak{m})=(t^4,t^7) \subseteq (t^4,t^5)$.  In fact,
for all $n \geq 4$,
$$(t^n)^*=(t^n,t^{n+1}) \supsetneq (t^n,t^{n+3})=(\frak{m}(t^n):\frak{m}).$$
\end{exm}

\end{document}